\documentclass[12pt]{article}

\usepackage{amsmath,amssymb,amsthm,amsfonts}
\usepackage{nicematrix,tikz}
\usetikzlibrary{fit}
\tikzset{highlight/.style={rectangle,
fill=gray!50,
rounded corners = 0.5 mm, 
inner sep=1pt,
fit=#1}}

\def\spn{\mathop{\rm span}}

\def\ind{\mathop{\rm ind}}
\def\ad{\mathop{\rm ad}}

\def\ind{{\rm ind}}

\def\phi{\varphi}
\def\g{\mathfrak g}
\def\i{\mathfrak i}
\def\m{\mathfrak m}

\def\F{\mathbb F}
\def\Z{\mathbb Z}
\def\N{\mathbb N}
\def\P{\mathbb P}

\makeatletter \let\@@pmod\pmod
\DeclareRobustCommand{\pmod}{\@ifstar\@pmods\@@pmod}
\def\@pmods#1{\mkern4mu({\operator@font mod}\mkern 6mu#1)}
\makeatother

\newtheorem{theorem}{Theorem}[section]
\newtheorem{lemma}[theorem]{Lemma}

\newtheorem{remark}{{\bf Remark}}
\providecommand{\keywords}[1]{\noindent{Keywords:} #1}
\providecommand{\classify}[1]{\noindent{Mathematics Subject
    Classification:} #1}

\title{Central Extensions of Restricted
  Affine Nilpotent Lie Algebras $n_+(A^{(1)}_1)(p)$}

\author{Tyler J. Evans \\Department of Mathematics \\California State
  Polytechnic University - Humboldt \\Arcata, CA 95521 USA
  \\evans@humboldt.edu \and Alice Fialowski \\Faculty of Informatics
  \\ E\" otv\" os Lor\' and University \\
  Budapest, Hungary\\ alice.fialowski@gmail.com}

\date{}

\begin{document}
\maketitle

\begin{center}
 {\it We dedicate this paper to Karl Hofmann \\ on the occasion of his 90th birthday.}
\end{center}

\bigskip
\begin{abstract}
  Consider the maximal nilpotent subalgebra $n_+(A_1^{(1)})$ of the
  simplest affine algebra $A_1^{(1)}$ which is one of the $\N$-graded
  Lie algebras with minimal number of generators. We show
  truncated versions of this algebra in positive characteristic admit the structure of
  a family of restricted Lie algebras. We compute the ordinary and
  restricted 1- and 2-cohomology spaces with trivial coefficients by
  giving bases. With these we explicitly describe the restricted
  1-dimensional central extensions.
\end{abstract}

{\footnotesize

\classify{17B50, 17B56,17B67, 17B70}}

\section{Introduction}

$\N$-graded Lie algebras over a field $\F$ are a special important
class of infinite dimensional Lie algebras. By $\N$-graded we mean
that the Lie algebra $\g$ is the direct sum of subspaces
$\g_i, i \in \N$, such that $[\g_i,\g_j] \subset \g_{i+j}$. Such
algebras obviously must have at least two generators. A special class
of these algebras are Lie algebras of maximal class, where
dim$ \g_1=2$, for $i \geq 2$, dim$\g_i=1$ and $[\g_1, \g_i]=\g_{i+1}$
for $i \geq 1$. These are also called filiform Lie algebras, and it is
known (see \cite{Fi}) that exactly three of them are $\N$-graded. If
we denote the basis elements by ${e_i, i \in \N}$, then these three
are the maximal nilpotent subalgebra $L_1$ of the Witt algebra with
nonzero brackets $[e_i,e_j]=(j-i)e_{i+j}$, the algebra $\m_0$ with
nonzero brackets $[e_1,e_i]=e_{i+1}$ for $i \geq 2$ and the Lie
algebra $\m_2$ with nonzero brackets $[e_1,e_i]=e_{i+1}$ for all
$i \geq 2$ and $[e_2,e_j]=e_{j+2}$ for all $j \geq 3$.

The filiform algebras $L_1, \m_0$ and $\m_2$ are well-studied, and
some of their important invariants are computed, such as the 1- and
2-cohomology spaces with trivial and adjoint coefficients. Finite
dimensional versions of these algebras are also studied in
characteristic $p > 0$. Namely, in \cite{ EvFim2, EvFim0, EvFiPe} we
considered restricted algebras of these types, and computed the
ordinary and restricted 1- and 2-cohomology spaces with trivial
coefficients.

There are other types of $\N$-graded Lie algebras with the minimal
number of generators $e_1$ and $e_2$. One of these is the maximal
nilpotent part $n_+(A^{(1)}_1)$ of the simplest affine Lie algebra
$A^{(1)}_1$ (see \cite{Kac}) with Cartan matrix
\[
  \begin{pmatrix}
    2 & -2\\
    -2 & 2
  \end{pmatrix}
\]
which has a basis $\{e_i\ |\ i\in\N\}$ with brackets
\begin{equation} [e_i,e_j]=\left\{
    \begin{array}{cl}
      e_{i+j} & \mbox{\rm if $j-i \equiv 1 \pmod* 3$;}\\
      0      & \mbox{\rm if $j-i \equiv 0 \pmod* 3$;}\\
      -e_{i+j} & \mbox{\rm if $j-i \equiv -1 \pmod* 3$.}
    \end{array}   
  \right.
\end{equation}    
 
The other $\N$-graded nilpotent affine Lie algebra is the maximal
nilpotent subalgebra of $n_+(BA_2)$ (see \cite{Kac}), which is defined
by the Cartan matrix
\[
  \begin{pmatrix}
    2 & -4\\
    -1 & 2
  \end{pmatrix}.
\]
  
Here there is a basis $\{e_i\ |\ i\in\N\}$ with brackets
\[
  [e_i,e_j]=b_{ij}e_{i+j}
\]
where the coefficients $b_{ij}$ depend on the residue obtained when
dividing $i$ and $j$ by 8 according to the rule
\[
  b_{ij}+b_{i'j'}=0 \ \mbox{\rm if $i+i'$ and $j+j'$ are divisible by
    $8$.}
\]
The Table below gives the numbers $b_{ij}$ (the others are determined
from $b_{ij}=-b_{ji}$ and $b_{ij}+b_{8-i,8-j}=0$).

\begin{center}
  \begin{NiceTabular}{|rr|ccccccc|}
    \hline
    \Block{1-*}{$i\pmod* 8$}\\
    &\Block[r]{1-2}{1} &  & 2 & 3 & 4 & 5 & 6 & 7\\
    \hline
    \Block[c]{4-1}{\rotate $j\pmod* 8$} &  0 & 1 & -2 & -1 & 0 & 1 & 2 & -1\\
    &1 & 0 & 1 & -1 & 3 & -2 & 0 & 1\\
    &2 & -1 & 0 & 0 & 0 & 1 & -1 & 0\\
    &3 & 1 & 0 & 0 & 3 & 1 & 1 & -2\\
    \hline
  \end{NiceTabular}
\end{center}

In addition to these five $\N$-graded algebras with two generators
there is a family of Lie algebras with countably many parameters
$\lambda_{4k} \in \F\P^1$, see \cite{Fi}.
 
The cohomology with trivial coefficients of the affine subalgebras
$n_+(A^{(1)}_1)$ and $n_+(BA_1)$ are computed, but nothing is known
about restricted cohomology, even though it gives an important
invariant of those algebras.
 
As in the case of the three filiform algebras with two generators, it
is interesting to consider whether there is a restricted Lie algebra
structure on the truncated algebras $n_+(A^{(1)}_1)(p)$ where $p>0$ is
a prime, and to describe their central extensions, as we did for the
three filiform cases \cite{EvFim2,EvFim0, EvFiPe}.
 
In this paper we concentrate on the algebras $n_+(A^{(1)}_1)(p)$. The
goal of the paper is to classify the restricted one-dimensional
central extensions of these restricted algebras. To do this we compute
the ordinary and restricted 2-cohomology spaces with trivial
coefficients. For completeness, we also compute the the ordinary and
restricted 1-cohomology spaces.
 
The structure of the paper is as follows. In Section 1 we define
finite dimensional truncated versions of the algebra $n_+(A^{(1)}_1)$
in positive characteristic $p>0$ and show that they admit the
structure of (a family of) restricted Lie algebras. Explicit formulas
for the Lie brackets and $[p]$-operators are given. In Section 2 we
compute the ordinary 1- and 2-cohomology spaces with trivial
coefficients giving descriptions of the bases.  Section 3 contains the
computation of the restricted 1- and 2-cohomology spaces with trivial
coefficients, again with descriptions of the bases. Finally, in
Section 4 we explicitly describe the restricted one-dimensional
central extensions corresponding to the non-trivial restricted
2-cocycles.

\section{Restricted Lie Algebra Structures}

Let $p$ be a prime integer and $\F$ a field of characteristic $p$. For $i,j\in\N$, define
\begin{equation}
  \label{eq:1}
  a_{i,j}=\left\{
    \begin{array}{cl}
      -1 & \mbox{\rm if $j-i\equiv -1\pmod* 3$;}\\
      0 & \mbox{\rm if $j-i\equiv 0\pmod* 3$;}\\
      1 & \mbox{\rm if $j-i\equiv 1\pmod* 3$.}
    \end{array}
  \right.
\end{equation}
Throughout the paper, we will let $[m]_3$ denote the congruence class
of an integer $m$ modulo $3$ so $a_{i,j}=[j-i]_3$. Let
$\g=n_+(A^{(1)}_1)=\bigoplus_{i=1}^\infty \F e_i$ with bracket
$[e_i,e_j]=a_{i,j}e_{i+j}$, $\i=(e_{p+1})$ the ideal generated by
$e_{p+1}$ and
\[\g (p)=n_+(A^{(1)}_1)(p)=n_+(A^{(1)}_1)/\i.\]
Then $\g (p)$ is a finite
dimensional, $\N$-graded Lie algebra with basis $\{e_1,\dots, e_p\}$
and $k$-th graded component $\g_k (p)=\F e_k$, $1\le k\le p$ and
$\g_k(p)=0$, $k>p$. If $\alpha_i,\beta_i\in\F$,
$g=\sum_{i=1}^p\alpha_ie_i$ and $h=\sum_{i=1}^p\beta_ie_i$, then
\begin{equation}
  \label{eq:2}
  [g,h]=\sum_{i=1}^{(p-1)/2}\sum_{j=i+1}^{p-i}
  a_{i,j}(\alpha_i\beta_j-\beta_i\alpha_j) e_{i+j}.
\end{equation}
For $g_1,\dots, g_n\in\g (p)$, denote the $n$-fold bracket by
\[[[\cdots[ [g_1,g_2],g_3]\cdots ],g_n]:=[g_1,g_2,\dots, g_n] .\] 
Note that $[g_1,g_2,\dots, g_p]=0$ for all $g_1,\dots, g_p\in\g (p)$
so for all $g,h\in\g (p)$, we have
\begin{equation*}
  (\ad g)^p(h) = [h,\underbrace{g,\dots, g}_p]=0.
\end{equation*}
If, for each $1\le k\le p$, we choose an element $e_k^{[p]}$ in the
center of $\g (p)$, we have
\begin{equation*}
  \ad e_k^{[p]} = 0 = (\ad e_k)^p,
\end{equation*} 
so the operator $-{}^{[p]}:\g (p)\to\g (p)$ gives $\g (p)$ the structure of a
restricted Lie algebra \cite{Jac, StradeFarnsteiner}. Since $p$-fold
brackets are $0$ in $\g (p)$, we have $(g+h)^{[p]} = g^{[p]}+h^{[p]}$ for all
$g,h\in\g (p)$. The center $Z(\g (p))$ is given by
\begin{equation*}
  Z(\g (p))=
  \left\{
    \begin{array}{cl}
      \F e_p & \mbox{\rm if $p\not\equiv 2\pmod* 3$};\\
      \F e_{p-1}\oplus \F e_p & \mbox{\rm if $p\equiv 2\pmod* 3$}.
    \end{array}
  \right.
\end{equation*}
If $p\not\equiv 2\pmod* 3$, then choosing $e_k^{[p]}\in Z(\g (p))$ is
equivalent to choosing $\lambda_k\in\F$ and setting
$e_k^{[p]}=\lambda_k e_p$ hence $\lambda\in\F^p$ determines a restricted
Lie algebra structure on $\g (p)$ which we denote by $\g
^\lambda(p)$. Conversely, if $-{}^{[p]}:\g (p)\to\g (p)$ is any
restricted Lie algebra structure on $\g (p)$, then
$\ad e_k^{[p]} = (\ad e_k)^p=0$ for all $1\le k\le p$ so that
$e_k^{[p]}=\lambda_k e_p$ for some $\lambda_k\in\F$.  For any such
restricted Lie algebra structure, if $g=\sum_{i=1}^p\alpha_ie_i$, then
\begin{equation}
  \label{eq:6}
  g ^{[p]} = \left (\sum_{i=1}^p \alpha_i^p\lambda_i\right) e_p.
\end{equation}
Likewise, if $p\equiv 2\pmod* 3$, then the restricted Lie algebra
structures on $\g(p)$ correspond to pairs of elements
$\mu,\lambda\in\F^p$ where $e_k^{[p]}=\mu_ke_{p-1}+\lambda_k e_p$. 
Denote this restricted Lie algebra by $\g ^{\mu,\lambda}(p)$. If
$g=\sum_{i=1}^p\alpha_ie_i$, then
\begin{equation}
  \label{eq:7}
  g ^{[p]} = \left (\sum_{i=1}^p \alpha_i^p\mu_i\right) e_{p-1}+\left (\sum_{i=1}^p \alpha_i^p\lambda_i\right) e_p.
\end{equation}

If $p=2$, then $n_+(A^{(1)}_1)(2)=\m_0(2)$, and if $p=3$,
$n_+(A^{(1)}_1)(3)=\m_0(3)=\m_2(3)=L_1(3)$ and hence the possible
restricted structures are the same.  These filiform algebras were
already studied in \cite{EvFim2,EvFim0} so everywhere below we assume
$p\ge 5$.

\section{Ordinary 1- and 2-Cohomology}

Our primary interest is in classifying (restricted) one-dimensional
central extensions, so we describe only the (graded) cochain spaces
$C^q=C^q(\g (p))=C^q(\g (p),\F)$ for $q=0,1,2,3$ and differentials
$d^q:C^q(\g (p))\to C^{q+1}(\g (p))$ for $q=0,1,2$ (for more details
on this cochain complex we refer the reader to
\cite{ChE,FuchsBook}). Set $C^0=\F$ and $C^q= (\wedge^q\g (p))^*$ for
$q=1,2,3$. We will use the following bases, ordered lexicographically,
throughout the paper.
\begin{align*}
  C^0:& \{1\}&\\
  C^1:&  \{e^k\ |\ 1\le k\le p\}&\\
  C^2:& \{e^{i,j}\ |\ 1\le i<j\le p\}&\\
  C^3:& \{e^{u,v,w}\ |\ 1\le u<v<w\le p\}&\\
\end{align*}
Here $e^k$, $e^{i,j}$ and $e^{u,v,w}$ denote the dual vectors of the
basis vectors $e_k\in \g(p)$,
$e_{i,j}=e_i\wedge e_j\in \wedge^2\g (p) $ and
$e_{u,v,w}=e_u\wedge e_v\wedge e_w\in \wedge^3\g (p)$, respectively.
The differentials $d^q:C^q\to C^{q+1}$ are defined for $\psi\in C^1$,
$\phi\in C^2$ and $g,h,f\in\g$ by\small
\begin{align*}
  d^0: C^0\to C^1,  &\  d^0=0&\\
  d^1:C^1\to C^2, &\   d^1(\psi)(g\wedge h)=\psi([g,h])&\\
  d^2:C^2\to C^3, &\  d^2(\phi)(g\wedge h\wedge f)=\phi([g,h]\wedge f)-\phi([g,f]\wedge h)+\phi([h,f]\wedge g). &\\
\end{align*}\normalsize
The cochain spaces $C^q(\g (p))$ are graded, and\small
\begin{align*}
  C^1_k(\g (p))&=\spn\{e^k\ |\ 1\le k\le p\}& \\
  C^2_k(\g (p))&=\spn\{e^{i,j}\ |\ 1\le i<j \le p, i+j=k, 3\le k\le 2p-1\}&  \\
  C^3_k(\g (p))&=\spn\{e^{u,v,w}\ |\ 1\le u<v<w \le p, u+v+w=k, 6\le k\le 3p-3\}&  
\end{align*}\normalsize
We have $\dim(C^1_k)=1$ for $1\le k\le p$, $\dim(C^2_k)=s(k)$ for
$3\le k\le p+1$ and $\dim(C^2_k)=\dim (C^2_{2p+2-k})$ for
$p+2\le k\le 2p-1$ where the map $s:\N\to\N$ is defined by
\begin{equation*}
  s(k)=\left\{
    \begin{array}{cl}
      \frac{k}{2}-1 & \mbox{\rm if $k$ is even;}\\
      \frac{k-1}{2} & \mbox{\rm if $k$ is odd.}\\
    \end{array}
  \right.
\end{equation*}
The differentials $d^q_k:C^q_k\to C^{q+1}_k$ are graded maps so
$H^q = \bigoplus_k H^q_k$ for $q=1,2$. We note that $d^1_1=d^1_2=0$
and $d^2_3=d^2_4=d^2_5=0$.  If $3\le k\le p$, direct computation shows
that
\begin{equation}
  \label{eq:9}
  d^1_k(e^k)=\sum_{i=1}^{s(k)} a_{i,k-i} e^{i,k-i}.
\end{equation}
The expression (\ref{eq:9}) together with $d^1_1=d^1_2=0$ immediately
gives the following
\begin{theorem}
  \label{ordinaryH1}
  A basis for $\ker d^1$ is $\{e^1,e^2\}$, hence $\dim H^1(\g (p))=2$,
  and the classes of $\{e^1,e^2\}$ form a basis.
\end{theorem}

The situation for $H^2_k$ is more complicated. If $6\le k\le 2p-1$,
$e^{i,j}\in C^2_k(\g (p))$ (so $i+j=k$), then the formula for the
differential $d^2$ and the restricted bracket (\ref{eq:2}) gives
\begin{align}
  \label{eq:10}
  \begin{split}
    d^2_k(e^{i,j}) =&\sum_{n=1}^{s(i)} a_{n,i-n} e^{n,i-n,k-i} -
    \sum_{n=k-2i+1}^{s(k-i)} a_{n,k-i-n} e^{n, k-i-n,
      i}\\
    +& \sum_{n=1}^{i-1} a_{n,k-i-n} e^{n,i,k-i-n} - \sum_{n=i+1}^{p-2}
    a_{n,k-i-n}
    e^{i,n,k-i-n}.\\
  \end{split}
\end{align}
We will use the following notation to describe index ranges in our
computations:
\begin{equation*}
  M(p,k)=\left\{
    \begin{array}{cl}
      1 & \mbox{\rm if $3\le k\le p+1$;}\\
      k-p & \mbox{\rm if $p+2\le k \le 2p-1$,}\\
    \end{array}
  \right.
\end{equation*}
\begin{equation*}
  G(p,k)=\left\{
    \begin{array}{cl}
      1 & \mbox{\rm if $6\le k\le 2p$;}\\
      1+ k-2p & \mbox{\rm if $2p+1\le k\le 3p-3$,}\\
    \end{array}
  \right.
\end{equation*}
\begin{equation*}
  F(p,k)=\left\{
    \begin{array}{cl}
      \frac{k}{3}-1 & \mbox{\rm if $k\equiv 0\pmod* 3$;}\\
      \frac{k-1}{3}-1 & \mbox{\rm if $k\equiv 1\pmod* 3$;}\\
      \frac{k-2}{3}-1 & \mbox{\rm if $k\equiv 2\pmod* 3$.}
    \end{array}
  \right.
\end{equation*}
If $e^{u,v,k-v-u}\in C^3_k$, then $G(p,k)\le u\le F(p,k)$ and
$u+1\le v\le s(k-u)$. If
\[d^2_k\left (\sum_{\substack{1\le i<j\le p\\i+j=k}} \alpha_{i,j}
    e^{ij}\right ) =\sum_{\substack{1\le i<j\le p\\i+j=k}}
  \alpha_{i,j}\sum_{\substack{1\le u< v< w\le p\\u+v+w=k}}
  A^{ij}_{uvw} e^{u,v,w}=0,\] then the term $e^{u,v,k-v-u}$ gives the
``{\it basic equation}''\small
\begin{equation}
  \label{system}
  \boxed{
    \underbrace{-\alpha_{u,k-u} a_{v,k-v-u}}_{k-u\le p} +
    \underbrace{\alpha_{v,k-v} a_{u,k-v-u}}_{k-v\le p} +
    \underbrace{\underbrace{\alpha_{u+v,k-v-u} a_{u,v}}_{u+v<k-v-u\le
        p}- \underbrace{\alpha_{k-v-u,u+v} a_{u,v}}_{k-v-u<u+v\le
        p}}_{\mbox{\rm\scriptsize never both non-zero}}=0.}
\end{equation}\normalsize
The kernel of $d^2_k$ is the space of solutions to the system
(\ref{system}) where $G(p,k)\le u\le F(p,k)$ and
$u+1\le v \le s(k-u)$. For $3\le k \le 2p-1$, define the cocycles
\begin{equation}
  \label{eq:15}
  \phi_k = \sum_{i=M(p,k)}^{s(k)} a_{i,k-i} e^{i,k-i}.
\end{equation}
The next Theorem is the main result of this section.
\begin{theorem}
  \label{ordinaryH2}
  A basis for $\ker d^2$ is
  \[\{e^{1,4}, e^{2,5}, \phi_3, \phi_4, \dots, \phi_{p+1}\}\] so
  $\dim\ker d^2 = p+1$. Moreover, $\phi_k=d^1(e^k)$ for $3\le k\le p$,
  hence $\dim H^2(\g (p))=3$ and the classes of
  $\{e^{1,4}, e^{2,5}, \phi_{p+1}\}$ form a basis.
\end{theorem}
We have already remarked that the cochain spaces and differentials are
graded, so Theorem~\ref{ordinaryH2} follows from 
\begin{theorem}
  \label{gradedk}
  If $k\ge p+2$, then $\ker d^2_k=0$. If $k< p+2$, then
  $\ker d^2_k = \{\phi_k\}$ unless $k=5$ or $k=7$. For all $p\ge 5$,
  $\ker d^2_5 = \{e^{1,4},\phi_5\}$, and $\ker d^2_7 = \{e^{2,5}\}$
  when $p=5$ and $\ker d^2_7 = \{e^{2,5},\phi_7\}$ when
  $p>5$. Moreover $\phi_k=d^1(e^k)$ for $3\le k\le p$, hence
  $H^2_k(\g (p))=0$ unless $k=5,7,p+1$. Bases for $H^2_5(\g (p))$,
  $H^2_7(\g (p))$ and $H^2_{p+1}(\g (p))$ consist of the classes of
  $\{e^{1,4}\}$, $\{e^{2,5}\}$ and $\{\phi_{p+1}\}$, respectively.
\end{theorem}

\begin{proof} 
  We have already noted that $d^2_3=d^2_4=d^2_5=0$, and hence
  $\{e^{1,2}\}$, $\{e^{1,3}\}$ and $\{e^{1,4},e^{1,5}\}$ form bases
  for $\ker d^2_k$ for $k=3,4,5$, respectively.  A direct computation
  using formula (\ref{eq:10}) for $d^2$ gives a basis for $\ker d^2_6$
  which is $\{\phi_6\}$, a basis for $\ker d^2_7$ is $\{e^{2,5}\}$
  when $p=5$, and a basis for $\ker d^2_7$ is $\{e^{2,5},\phi_7\}$
  when $p>5$.

  For the remainder of the proof we assume $p\ge 5$ and $k\ge 8$. The
  proof is combinatorial and consists of three parts:

  \begin{itemize}
  \item [(i)] $k > p+2$;
  \item [(ii)] $k=p+2$;
  \item [(iii)] $k< p+2$.
  \end{itemize}

  In each part, we will treat the cases $k=6t+r$ where $0<t$ and
  $0\le r< 6$ separately, and only use the {\it basic equations}
  (\ref{system}) with $u=1$ or $u=2$. We have included examples for
  $p=23$ in the Appendix to illustrate the computations.

  \paragraph{(i)} $k > p+2\Longrightarrow \ker d^2_k=0$.

  If $u=1$ and $v=i-1$, then (\ref{system}) reduces to
  \begin{align}
    \begin{split}
      \label{bigksystemu1}
      &\alpha_{k-p,p}[r-p+1]_3 =0 \hspace{4.3cm} ( i=k-p)\\
      &\alpha_{i-1,k-i+1}[r-i-1]_3+\alpha_{i,k-i}[i+1]_3=0
      \hspace{1cm} (k-p+1\le i\le s(k)),
    \end{split}
  \end{align}
  and if $u=2$ and $v=i-2$, then (\ref{system}) reduces to
  \begin{align}
    \begin{split}
      \label{bigksystemu2}
      &\alpha_{k-p,p}[r-p-1]_3 =0 \hspace{4.3cm} ( i=k-p)\\
      &\alpha_{k-p+1,p-1}[r-p]_3 =0 \hspace{4.3cm} ( i=k-p+1)\\
      &\alpha_{i-2,k-i+2}[r-i+1]_3+\alpha_{i,k-i}[i-1]_3=0
      \hspace{1cm} (k-p+2\le i\le s(k)).
    \end{split}
  \end{align}
  Suppose $[r]_3=1$. If $[i]_3=0$, then (\ref{bigksystemu1}) shows
  $\alpha_{i,k-i}=0$. If $[i]_3=1$, then $[i-1]_3=0$ so that
  (\ref{bigksystemu1}) again shows $\alpha_{i,k-i}=0$. If $[i]_3=-1$,
  then (\ref{bigksystemu2}) shows $\alpha_{i,k-i}=0$.

  Suppose that $[r]_3=-1$.  If $[i]_3=1$, then (\ref{bigksystemu1})
  shows $\alpha_{i,k-i}=0$. If $[i]_3=0$, then (\ref{bigksystemu1})
  shows $\alpha_{i-1,k-i+1}+\alpha_{i,k-i}=0$ and (\ref{bigksystemu2})
  shows $\alpha_{i,k-i}=0$, hence $\alpha_{i-1,k-i+1}=0$ as well.  (If
  $r=5$, then $[s(k)]_3=-1\ne 0$, but $[s(k)-2]_3=0$ so again
  (\ref{bigksystemu2}) implies $\alpha_{s(k),k-s(k)}=0$.)

  Finally suppose $[r]_3=0$. If $[i]_3=0$, then (\ref{bigksystemu1})
  applied to $i$ and $i+1$ gives
  \begin{equation}
    \label{groupedbythree}
    \alpha_{i-1,k-i+1}=\alpha_{i,k-i}=\alpha_{i+1,k-i-1}.
  \end{equation}
  Moreover, using (\ref{bigksystemu2}) and (\ref{groupedbythree}), we
  see that if $k-p+2\le i$ and $\alpha_{i-2,k-i+2}=0$, then
  $\alpha_{i-1,k-i+1}=\alpha_{i,k-i}=\alpha_{i+1,k-i-1}=0$. If
  $[p]_3=-1$, then (\ref{bigksystemu1}) shows $\alpha_{k-p,p}=0$ and
  (\ref{bigksystemu2}) shows $\alpha_{k-p+1,p-1}=0$, and hence
  $\alpha_{i,k-i}=0$ for all $i$ by the previous remark. If $[p]_3=1$,
  then (\ref{bigksystemu2}) shows $\alpha_{k-p,p}=0$ so again,
  $\alpha_{i,k-i}=0$ for all $i$.  This completes the proof of part
  (i).

  \paragraph{(ii)} $k = p+2 \Longrightarrow \ker d^2_k=0$

  Since $p\ge 5$ is a prime, $p+2=k=6t+r$ implies $r=1$ or
  $r=3$. Moreover, $r=1$ if and only if $[p]_3=-1$ and $r=3$ if and
  only if $[p]_3=1$. Suppose that $r=1$. If $u=1$ and $v=i-1$, then
  the {\it basic equation} (\ref{system}) still reduces to
  (\ref{bigksystemu1}). If $u=2$ and $v=i-2$, then (\ref{system})
  reduces to
  \begin{align}
    \begin{split}
      \label{p+2systemu2}
      &-\alpha_{2,p}[i]_3+\alpha_{i-2,k-i+2}[r-i+1]_3+\alpha_{i,k-i}[i-1]_3=0\\
      & \hspace{2.5cm} (k-p+3\le i\le s(k))
    \end{split}
  \end{align}
  In particular, taking $i=5$ and $i=7$, we have, respectively,
  \[\alpha_{2,p}+\alpha_{5,p-3}=0\]
  and
  \[-\alpha_{2,p}+\alpha_{5,p-3}=0.\] It follows that $\alpha_{2,p}=0$
  and hence (\ref{p+2systemu2}) reduces to (\ref{bigksystemu2}), and
  the argument for $[r]_3=1$ in part (i) shows that $\alpha_{i,k-i}=0$ for
  all $i$.

  Suppose $r=3$.  If $u=1$ and $v=i-1$, then the {\it basic
    equation} (\ref{system}) again reduces to (\ref{bigksystemu1}),
  and applying (\ref{bigksystemu1}) to $i$ and $i+1$ for $[i]_3=0$, we
  have
  \begin{equation}
    \label{groupby3p+2}
    \alpha_{i-1,k-i+1}=\alpha_{i,k-i}=\alpha_{i+1,k-i-1}.
  \end{equation}
  Let $A_m=\alpha_{3m,k-3m}$ for $1\le m\le t$. Now, we use
  (\ref{p+2systemu2}) for all $[i]_3=1$ to get a system of equations
  \begin{align}
    \begin{split}
      \label{auxsystem}
      -A_1-A_1+A_2 &= 0\\
      -A_1-A_m+A_{m+1}&= 0 \hspace{.5cm} (2\le m\le t-1)\\
      -A_1-A_t-A_t &=0.
    \end{split}
  \end{align}
  These equations imply $A_t=-tA_1$ and $(2t-1)A_1=0$ so that
  $A_1=A_t=0$. From this it follows that $A_m=0$ for all $m$ and hence
  (\ref{groupby3p+2}) gives
  $\alpha_{i,k-i}=0$ for all $i$. This completes the proof of part
  (ii).

  \paragraph{(iii)}
  $k < p+2\Longrightarrow \ker d^2_k=\spn_\F\{\phi_k\}$

  Recall $\dim C^2_k(\g(p)) = s(k)$ for $k<p+2$. A direct calculation
  shows that if $k< p+2$, then the element $\phi_k$ defined in
  (\ref{eq:15}) is a solution of the {\it basic equation}
  (\ref{system}) so that $\phi_k\in\ker(d^2_k)$. We complete the proof
  by showing that the set
  \begin{equation}
    \label{theset}
    \{d^2_k(e^{i,k-i})\ |\ 1\le i\le s(k)-1\}
  \end{equation}
  is linearly independent.

  If $u=1$ and $v=i$ satisfies $2\le i\le s(k)-1$, then (\ref{system})
  reduces to
  \begin{equation}
    \label{smallksystemu1}
    -\alpha_{1,k-1}[r+i-1]_3 +
    \alpha_{i,k-i}[r-i+1]_3+\alpha_{i+1,k-i-1}[i-1]_3=0.
  \end{equation}
  If $u=2$ and $v=i$ satisfies $3\le i\le s(k)-2$, then (\ref{system})
  reduces to
  \begin{equation}
    \label{smallksystemu2}
    -\alpha_{2,k-2}[r+i+1]_3 +
    \alpha_{i,k-i}[r-i-1]_3+\alpha_{i+2,k-i-2}[i-2]_3=0.
  \end{equation}
  We begin by showing that if $[r]_3=1$ or $2$, then
  $\alpha_{1,k-1}=\alpha_{2,k-2}=0$ and hence (\ref{smallksystemu1})
  and (\ref{smallksystemu2}) reduce to (\ref{bigksystemu1}) and
  (\ref{bigksystemu2}), respectively. The proof in (i) then shows
  $\alpha_{i,k-i}=0$ for $2\le i\le s(k)-1$ and the set (\ref{theset})
  is linearly independent.

  The following Table lists selected values of $v$ for each value of
  $r$ and the resulting simplified equations from
  (\ref{system}). Since we are considering only linear combinations of
  the elements of the set (\ref{theset}), these equations show
  $\alpha_{1,k-1}=\alpha_{2,k-2}=0$ by inspection.

\[\begin{array}{|l|l|r|}
    \hline
    \multicolumn{3}{|c|}{u=1} \\
    \hline
    r=1 & v= s(k)-1 & \alpha_{1,k-1} +\alpha_{s(k),k-s(k)}=0\\
    \hline
    r=2,4 & v= s(k) & -\alpha_{1,k-1} +\alpha_{s(k),k-s(k)}[r/2-1]_3=0\\
    \hline
    r=5 & v=s(k)-2 &  -\alpha_{1,k-1} -\alpha_{s(k)-1,k-s(k)+1}=0\\
        & v=s(k)-1 & \alpha_{1,k-1} -\alpha_{s(k)-1,k-s(k)+1}=0 \\
    \hline
    \multicolumn{3}{|c|}{u=2} \\
    \hline
    r=1 & v=s(k)-3 &  \alpha_{2,k-2}+\alpha_{s(k)-1,k-s(k)+1}=0\\
        & v=s(k)-1 & -\alpha_{2,k-2}+\alpha_{s(k)-1,k-s(k)+1}=0\\
    \hline
    r=2 & v=s(k)-2 & -\alpha_{2,k-2}+\alpha_{s(k),k-s(k)}=0\\
    \hline
    r=4 & v=s(k)-1 & -\alpha_{2,k-2} = 0\\
    \hline
    r=5 &  v=s(k)-1 & -\alpha_{2,k-2}+\alpha_{s(k),k-s(k)}=0\\
    \hline
  \end{array}
\]

If $r=0$, then selecting $u=1$ and $v=s(k)$ shows $\alpha_{1,k-1}=0$
as above. Using this, (\ref{smallksystemu1}) reduces to
(\ref{bigksystemu1}), and applying (\ref{bigksystemu1}) to $i$ and
$i+1$ for $[i]_3=0$, we again have
\begin{equation}
  \label{groupby3}
  \alpha_{i-1,k-i+1}=\alpha_{i,k-i}=\alpha_{i+1,k-i-1}.
\end{equation}
Let $A_m=\alpha_{3m,k-3m}$ for $1\le m\le t-1$. Letting $u=2$ and
$v=i$ ($3\le i\le s(k)-2$) in the {\it basic equation}
(\ref{system}), we have
\begin{equation}
  \label{systemu2smallp}
  -\alpha_{2,k-2}[i+1]_3-\alpha_{i,k-i}[i+1]_3+\alpha_{i+2,k-i-2}[i+1]_3=0.
\end{equation}
Using (\ref{systemu2smallp}) for all $[i]_3=0$, we get a system of
equations
\begin{align*}
  \begin{split}
    -A_1-A_m+A_{m+1}&= 0 \hspace{.5cm} (1\le m\le t-2)\\
    -A_1-A_{t-1} &=0.
  \end{split}
\end{align*}
These equations imply $tA_1=0$ and hence $A_m=0$ for all $m$. This
together with (\ref{groupby3}) implies $\alpha_{i,k-i}=0$ for all $1\le i\le s(k)-1$.

Finally, we assume $r=3$. Letting $u=1$ and $u=2$, respectively, with
$v=i=s(k)-1$ in (\ref{system}), we have
\begin{align*}
  \alpha_{1,k-1}+\alpha_{s(k)-1,k-s(k)+1}-\alpha_{s(k),k-s(k)} & = 0\\
  -\alpha_{2,k-2}-\alpha_{s(k)-1,k-s(k)+1}-\alpha_{s(k),k-s(k)} & = 0 
\end{align*}
We are considering linear combinations of elements of (\ref{theset})
hence these equations show that
\begin{align}
  \label{theinductionseed}
  \alpha_{1,k-1}=-\alpha_{s(k)-1,k-s(k)+1}=\alpha_{2,k-2}.
\end{align}
Now, for $u=1$ and $u=2$ respectively, we have
\begin{align}
  \label{induction}
  \begin{split}
    \alpha_{i+1,k-i-1}[i-1]_3 & =
    \alpha_{1,k-1}[i-1]_3+\alpha_{i,k-i}[i-1]_3\ (2\le i\le s(k)-1)\\
    \alpha_{i+2,k-i-}[i+1]_3 & =
    \alpha_{2,k-2}[i+1]_3+\alpha_{i,k-i}[i+1]_3 \ (3\le i\le s(k)-2)
  \end{split}
\end{align}
Together, (\ref{theinductionseed}), (\ref{induction}) and an induction
shows that
\[\alpha_{i,k-i}=c(i)\alpha_{1,k-1}, 2\le i\le s(k)-1.\]
Moreover, $c(i)\le i$ for all $i$, hence, in particular,
\[c(s(k)-1)\le s(k)-1=3t<6t+1\le p-1.\] This implies that
\[\alpha_{s(k)-1,k-s(k)+1}=-\alpha_{1,k-1}\ \mbox{\rm and}\ 
\alpha_{s(k)-1,k-s(k)+1}=c(s(k)-1)\alpha_{1,k-1}\] with
$c(s(k)-1)\ne -1$ in $\F$. This is a contradiction unless
$\alpha_{1,k-1}=0$, so (\ref{theset}) is linear independent and the
proof is complete.
\end{proof}

\begin{remark}
  It is interesting to note that $\dim H^2(n_+(A_1^{(1)})(p))=3$ is
  independent of the prime $p$, just as for the filiform algebras
  $\m_2(p)$ \cite{EvFim2}. On the other hand, $\dim H^2 (\m_0(p))$
  depends on $p$ \cite{EvFim0}. It would be interesting to see whether
  or not the dimension of the 2-cohomology for the other affine
  algebras $n_+(BA_2)(p)$ depends on $p$. Of course, in all cases, the
  dimension of the ordinary 1-cohomology spaces is 2 as all of these
  algebras are 2-generated.
\end{remark}

\section{Restricted Cohomology}
For general information on restricted cohomology of restricted Lie
algebras, we refer the reader to \cite{Feldross,Jantzen2006,Jantzen}.
Where no confusion can arise, for notational simplicity we denote by
$C^q_*$ both the spaces
$C^q_*(\g^{\lambda} (p)) = C^q_*(\g^{\lambda} (p);\F)$ and
$C^q_*(\g^{\mu,\lambda} (p)) = C^q_*(\g^{\mu,\lambda} (p);\F)$. We
describe only the restricted cochain spaces $C^q_*$ for $q=0,1,2,3$
and the restricted differentials $d^q_*: C^q_*\to C^{q+1}_*$ for
$q=0,1,2$, and refer the reader to \cite{EvansFuchs2} for a more
detailed description of this partial cochain complex. The computations
of the restricted cohomology spaces $H^q_*(\g^{\lambda} (p);\F)$ and
$H^q_*(\g^{\mu,\lambda} (p);\F)$ for $q=1,2$ are carried out in
separate subsections.

Given $\phi\in C^2$, a map $\omega:\g (p)\to\F$ is {\bf
  $\phi$-compatible} (in \cite{EvansFuchs2}, the authors use the term
$\omega$ {\it has the $*$-property with respect to} $\phi$ ) if for
all $g,h\in\g (p)$ and all $\alpha\in\F$,
$\omega(\alpha g)=\alpha^p \omega (g)$ and
\begin{equation}
  \label{eq:16}
  \omega(g+h)=\omega(g)+\omega(h) + \sum_{\substack{g_i=\mbox{\rm\scriptsize $g$
        or $h$}\\ g_1=g, g_2=h}}
  \frac{1}{\#(g)}\phi([g_1,g_2,g_3,\dots,g_{p-1}]\wedge g_p).
\end{equation}
Note that $\omega:\g (p)\to\F$ is $0$-compatible if and only if
$\omega$ is $p$-semilinear.

If $\zeta\in C^3$, then a map $\eta:\g (p)\times \g (p)\to\F$ is {\bf
  $\zeta$-compatible} (in \cite{EvansFuchs2}, $\eta$ {\it has the
  $**$-property with respect to} $\zeta$ ) if for all $\alpha\in\F$
and all $g,h,h_1,h_2\in\g (p)$ we have $\eta(\cdot,h)$ linear in the
first coordinate, $\eta(g,\alpha h)=\alpha^p\eta(g,h)$ and
\begin{align*}
  \eta(g,h_1+h_2) &=
                    \eta(g,h_1)+\eta(g,h_2)-\nonumber \\
                  & \sum_{\substack{l_1,\dots,l_p=1 {\rm or} 2\\ l_1=1,
  l_2=2}}\frac{1}{\#\{l_i=1\}}\zeta (g\wedge
  [h_{l_1},\cdots,h_{l_{p-1}}]\wedge h_{l_{p}}).
\end{align*}
Define the restricted cochain spaces as $C^0_*=C^0$, $C^1_*=C^1$,
\[C^2_*=\{(\phi,\omega)\ |\ \phi\in C^2, \omega:\g\to\F\ \mbox{\rm is
    $\phi$-compatible}\}\]
\[C^3_*=\{(\zeta,\eta)\ |\ \zeta\in C^3, \eta:\g\times\g\to\F\
  \mbox{\rm is $\zeta$-compatible}\}.\] If $\phi\in C^2$, we can assign
values $\omega(e_k)$ arbitrarily to elements $e_k$ of a basis
$\{e_1,\dots, e_p\}$ for $\g (p)$, set
$\omega(\alpha e_k)=\alpha^p\omega(e_k)$ for all $\alpha\in\F$ and use
(\ref{eq:16}) to determine a unique $\phi$-compatible map
$\omega:\g (p)\to\F$ \cite{EvFim0}. In particular, we can define
$\tilde\phi(e_k)=0$ for all $k$ and use (\ref{eq:16}) to determine a
unique $\phi$-compatible map $\widetilde\phi:\g (p)\to\F$. Note that, in
general, $\widetilde\phi\ne 0$ but $\widetilde\phi (0)=0$.

\begin{lemma}
  If $\phi_1,\phi_2\in C^2$ and $\alpha\in\F$, then
  $\widetilde{(\alpha\phi_1+\phi_2)} = \alpha\widetilde\phi_1 +
  \widetilde\phi_2$.
\end{lemma}

\begin{proof}
  An easy computation shows
  $(\alpha\widetilde\phi_1 + \widetilde\phi_2)(e_k) = 0$ for all $k$
  and $\alpha\widetilde\phi_1 + \widetilde\phi_2$ satisfies the
  compatibility condition (\ref{eq:16}) with
  $\phi=\alpha \phi_1 +\phi_2$. The result now follows from
  uniqueness.
\end{proof}

For $1\le k\le p$, define $\overline e^k:\g (p)\to\F$ by
\[\overline e^k \left(\sum_{n=1}^p \alpha_n e_n\right ) =
  \alpha_k^p.\] It is shown in \cite{EvansFuchs2} that
\[\dim C^2_* = \binom{p+2-1}{2}=\binom{p+1}{2}=\binom{p}{2}+p,\]
and hence it follows that
\begin{equation*}
  \{(e^{i,j},\widetilde{e^{i,j}})\ |\ 1\le i<j\le p\} \cup \{(0,\overline e^k)\ |\ 1\le k\le p\}
\end{equation*}
is a basis for $C^2_*$. We use this basis in all computations that
follow. 

\begin{remark}
  \label{tildemaps}
  The $(p-1)$-fold bracket in (\ref{eq:16}) always gives a multiple of
  $e_p$. For $p\ge 5$, if $k<p+1$, $\phi_k$ is identically zero on
  $e_p\wedge \g^\lambda(p)$ (or $e_p\wedge\g^{\mu,\lambda}(p)$), and
  hence $\widetilde{\phi_k}=0$ because $\widetilde{\phi_k}(e_i)=0$ for
  all $i$. Likewise $\widetilde{e^{1,4}}=0$ and
  $\widetilde{e^{2,5}}=0$, unless $p=5$. The restriction of
  $\phi_{p+1}$ to $e_p\wedge \g^\lambda(p)$ (or
  $e_p\wedge \g^{\mu,\lambda}(p)$) is equal to $e^{1,p}$, hence
  $\widetilde{\phi_{p+1}}=\widetilde{e^{1,p}}$. We can then compute,
  using (\ref{eq:16}):
  \begin{align*}
    \widetilde{e^{2,5}}\left(\sum_{i=1}^p\alpha_ie_i\right) & =
                                                              \frac{1}{2}\alpha_1^3\alpha_2^2\ 
                                                            \,\,  (p\ne 5),\\
    \widetilde{\phi_{p+1}}\left(\sum_{i=1}^p\alpha_ie_i\right) & =
                                                                 \widetilde{e^{1,p}}\left(\sum_{i=1}^p\alpha_ie_i\right)  =
                                                                 \alpha_1^{p-1}\alpha_2. 
  \end{align*}
\end{remark}

For $\psi\in C^1_*$, define the map $\ind^1(\psi):\g (p)\to\F$
by
\[\ind^1(\psi)(g)=\psi(g^{[p]}).\] The map $\ind^1(\psi)$ is
$d^1(\psi)$-compatible for all $\psi\in C^1_*$, and the differential
$d^1_*:C^1_*\to C^2_*$ is defined by
$d^1_*(\psi) = (d^1(\psi),\ind^1(\psi))$. For
$(\phi,\omega)\in C^2_*$, define the map
$\ind^2(\phi,\omega):\g (p)\times\g (p)\to\F$ by the
formula
\[\ind^2(\phi,\omega)(g,h)=\phi(g\wedge h^{[p]}).\] The map
$\ind^2(\phi,\omega)$ is $d^2(\phi)$-compatible for all $\phi\in C^2$,
and the differential $d^2_*:C^2_*\to C^3_*$ is defined by
$d^2_*(\phi,\omega) = (d^2(\phi),\ind^2(\phi,\omega))$ (see
\cite{EvansFuchs2} for details). We note that (with trivial
coefficients) if $\omega_1$ and $\omega_2$ are both $\phi$-compatible,
then $\ind^2(\phi, \omega_1)=\ind^2(\phi, \omega_2)$.

\begin{lemma}
  \label{swap}
  If $(\phi,\omega)\in C^2_*$ and $\phi=d^1(\psi)$ with
  $\psi\in C^1$, then $(\phi,\ind^1(\psi))\in C^2_*$ and
  $\ind^2(\phi,\omega)=\ind^2(\phi,\ind^1(\psi))$.
\end{lemma}

\begin{proof}
  We know that $\ind^1(\psi)$ is $d^1(\psi)$-compatible for all
  $\psi\in C^1=C^1_*$ \cite{EvansFuchs2}. If $\phi=d^1(\psi)$, then
  $(\phi,\ind^1(\psi))=(d^1(\psi),\ind^1(\psi))\in C^2_*$, and
  $\ind^2$ depends only on $\phi$ by the last sentence in the previous
  paragraph.
\end{proof}

We will use the following notation everywhere in the next two
subsections: let $g=\sum \alpha_i e_i$, $h=\sum \beta_i e_i$,
$\psi=\sum \gamma_i e^i$ and $\phi=\sum \sigma_{i,j}e^{i,j}$.

\subsection{Restricted Cohomology for $p\not\equiv 2\pmod* 3$}
Assume $p\ge 5$ is a prime and $p\not\equiv 2\pmod* 3$. Using the
$[p]$-operator (\ref{eq:6}) and the definitions of the maps $\ind^1$
and $\ind^2$, we have
\begin{equation}
  \label{eq:19}
  \ind^1(\psi)(g)=
  \gamma_p\left (\sum_{i=1}^p \alpha_i^p\lambda_i\right )
\end{equation}
and
\begin{equation}
  \label{eq:21}
  \ind^2(\phi,\omega)(g,h)=
  \left (\sum_{i=1}^p \beta_i^p\lambda_i\right
  )\left(\sum_{i=1}^{p-1}\sigma_{i,p}\alpha_i\right).
\end{equation}

\begin{theorem}
  \label{restrictedH1}
  $H^1_*(\g^\lambda(p))=H^1(\g (p))$, and the classes of the cocycles
  $\{e^1,e^2\}$ form a basis.
\end{theorem}

\begin{proof}
  Recall that $H^1_*$ is the subspace of $H^1$ consisting of those
  cocycles $\psi\in C^1$ for which $\ind^1(\psi)=0$ (see \cite{Ho} or
  \cite{EvansFuchs2}). By Theorem~\ref{ordinaryH1}, $\psi\in C^1$ is a
  cocycle if and only if $\gamma_k=0$ for $3\le k\le p$. Since
  $p\ge 5$, then (\ref{eq:19}) shows that $\ind^1(\psi)=0$ for all
  cocycles $\psi$.
\end{proof}

\begin{theorem}
  \label{H2forpnot2}
  The set
  \begin{equation}
    \label{kernelbasispnot2}
    \{(e^{1,4},\widetilde{e^{1,4}}), (e^{2,5},\widetilde{e^{2,5}}),
    (\phi_3,\widetilde{\phi_3}), \dots,
    (\phi_{p+1},\widetilde{\phi_{p+1}}),(0,\overline e^1),\dots,
    (0,\overline e^p)\} 
  \end{equation}
  forms a basis for the kernel of $d^2_*$. Moreover,
  $\dim H^2_*(\g^\lambda(p))=p+3$ and the classes of
  \[\{(e^{1,4},\widetilde{e^{1,4}}), (e^{2,5},\widetilde{e^{2,5}}),
    (\phi_{p+1},\widetilde{\phi_{p+1}}),(0,\overline e^1),\dots,
    (0,\overline e^p)\}\] form a basis.
\end{theorem}

\begin{proof}
  If $\lambda=0$, then (\ref{eq:21}) shows that
  $\ind^2(\phi,\omega)=0$ for all $(\phi,\omega)\in C^2_*$. Therefore,
  if $\phi\in C^2$ is an ordinary cocycle, $d^2_*(\phi,\omega)=(0,0)$
  for any $\phi$-compatible map $\omega$. Moreover, if
  $d^2_*(\phi,\omega)=(0,0)$, then $\phi$ must be an ordinary
  cocycle. This, together with our result on the ordinary cohomology
  (Theorem~\ref{ordinaryH2}) shows that (\ref{kernelbasispnot2}) is a
  basis for the kernel of $d^2_*$. If $3\le k\le p$, then
  $\phi_k=d^1(e^k)$, and Lemma~\ref{swap} shows that we can replace
  $\widetilde{\phi_k}$ with $\ind^1(e^k)$, so
  $d^1_*(e^k)=(\phi_k,\ind^1(e^k))$. This shows that
  $\dim H^2_*(\g^{\lambda} (p))=p+3$ and the classes of
  \[\{(e^{1,4},\widetilde{e^{1,4}}), (e^{2,5},\widetilde{e^{2,5}}),
    (\phi_{p+1},\widetilde{\phi_{p+1}}),(0,\overline e^1),\dots,
    (0,\overline e^p)\}\] form a basis.

  If $\lambda\ne 0$, then $\lambda_n\ne 0$ for some $1\le n\le p$.  If
  $(\phi,\omega)\in C^2_*$, then (\ref{eq:21}) shows
  \[\ind^2(\phi,\omega)(e_m,e_n)=\lambda_n\sigma_{m,p}\]
  for all $1\le m\le p-1$. It follows that $d^2_*(\phi,\omega)=(0,0)$
  if and only if $d^2(\phi)=0$ and
  $\sigma_{1,p}=\sigma_{2,p}=\cdots = \sigma_{p-1,p}=0$. By
  Theorem~\ref{ordinaryH2}, if $d^2(\phi)=0$, then $i+j\le p+1$ for
  all $i,j$ so $\sigma_{2,p}=\cdots = \sigma_{p-1,p}=0$. Moreover, the
  coefficient of $e^{1,p}$ in $\phi_{p+1}$ is
  $\sigma_{1,p}=a_{1,p}=p-1\pmod* 3= 0\pmod* 3$. This shows that
  $d^2_*(\phi,\omega)=(0,0)$ if and only if $d^2(\phi)=0$ and hence
  (\ref{kernelbasispnot2}) is a basis for $\ker d^2_*$ in this case as
  well. The proof now proceeds as in the case $\lambda = 0$.
\end{proof}

\subsection{Restricted Cohomology for $p\equiv 2\pmod* 3$}

Using formula (\ref{eq:7}) for the $p$-operator and the definitions of the maps $\ind^1$ and
$\ind^2$, we have
\begin{equation}
  \label{eq:20}
  \ind^1(\psi)(g)=  \gamma_{p-1}\left (\sum_{i=1}^p \alpha_i^p\mu_i\right )+\gamma_p\left (\sum_{i=1}^p \alpha_i^p\lambda_i\right )
\end{equation}
and
\begin{align}
  \begin{split}
    \label{eq:22}
    \ind^2(\phi,\omega)(g,h)= &\left (\sum_{i=1}^p
      \beta_i^p\mu_i\right
    )\left(\sum_{i=1}^{p-2}\sigma_{i,p}\alpha_i-\sigma_{p-1,p}\alpha_p\right)\\
    &+\left (\sum_{i=1}^p \beta_i^p\lambda_i\right
    )\left(\sum_{i=1}^{p-1}\sigma_{i,p}\alpha_i\right).
  \end{split}
\end{align}

The congruence class of $p$ modulo $3$ has no affect on the
$1$-cohomology, hence the proof of Theorem~\ref{restrictedH1} also shows

\begin{theorem}
  $H^1_*(\g^{\mu,\lambda}(p))=H^1(\g (p))$, and the classes of the
  cocycles $\{e^1,e^2\}$ form a basis.
\end{theorem}

\begin{theorem}
  \label{H2forpcong2mupluslambda0}
  If $\mu+\lambda=0$, then
  $H^2_*(\g^{\mu,\lambda}(p))=H^2_*(\g^\lambda(p))$.
\end{theorem}

\begin{proof}
  If $(\phi,\omega)\in C^2_*$, then (\ref{eq:22}) shows for
  $1\le n\le p$
  \begin{equation}
    \label{eq:24}
    \ind(\phi,\omega) (e_m,e_n)=\left\{
      \begin{array}{ll}
        (\mu_n+\lambda_n)\sigma_{m,p} & \mbox{\rm if $1\le m\le p-2$}\\
        (-\mu_n+\lambda_n)\sigma_{m,p} & \mbox{\rm if $m=p-1$}\\
        0 & m=p
      \end{array}
    \right.
  \end{equation}
  If $\mu+\lambda=0$, it follows that $d^2_*(\phi,\omega)=(0,0)$
  if and only if $d^2(\phi)=0$ and $\sigma_{p-1,p}=0$. But
  $d^2(\phi)=0$ implies $\sigma_{p-1,p}=0$ by Theorem~\ref{ordinaryH2}
  so $d^2_*(\phi,\omega)=(0,0)$ if and only if
  $d^2(\phi)=0$ and (\ref{kernelbasispnot2}) is a basis for
  $\ker d^2_*$. The argument now procedes precisely as in the proof of
  Theorem~\ref{H2forpnot2}.
\end{proof}

\begin{theorem}
  \label{H2forpcong2mupluslambdanot0}
  If $\mu+\lambda\ne 0$, then $\dim H^2_*(\g^{\mu,\lambda}(p))=p+2$
  and the classes of
  \[\{(e^{1,4},\widetilde{e^{1,4}}), (e^{2,5},\widetilde{e^{2,5}})
    ,(0,\overline e^1),\dots, (0,\overline e^p)\}\] form a basis.
\end{theorem}

\begin{proof}
  If $\mu+\lambda\ne 0$, then looking at (\ref{eq:24}), we see
  $d^2_*(\phi,\omega)=(0,0)$ implies $d^2(\phi)=0$ and
  $\sigma_{1,p}=\cdots = \sigma_{p-2,p}=0$. Theorem~\ref{ordinaryH2}
  shows that if $d^2(\phi)=0$, then
  $\sigma_{2,p}=\cdots =\sigma_{p-1,p}=0$. Moreover for $\phi_{p+1}$,
  $\sigma_{1,p}=a_{1,p}=p-1\not\equiv 0\pmod* 3$, and hence
  $d^2_*(\phi_{p+1},\widetilde{\phi_{p+1}})\ne (0,0)$. For any other
  ordinary cocycle $\phi$, we have $\sigma_{1,p}=0$. Therefore the kernel of
  $d^2_*$ is
  \begin{equation*}
    \{(e^{1,4},\widetilde{e^{1,4}}), (e^{2,5},\widetilde{e^{2,5}}),
    (\phi_3,\widetilde{\phi_3}), \dots,
    (\phi_{p},\widetilde{\phi_{p}}),(0,\overline e^1),\dots,
    (0,\overline e^p)\},
  \end{equation*}
  and we again procede as in the proof of
  Theorem~\ref{H2forpnot2}.
\end{proof}

\section{One-Dimensional Central Extensions}
One-dimensional central extensions $E=\g\oplus\F c$ of an ordinary Lie
algebra $\g$ are parameterized by the cohomology space $H^2(\g)$
\cite[Ch.~1, Sec.~4.6]{FuchsBook}, and restricted one-dimensional
central extensions of a restricted Lie algebra $\g$ with $c^{[p]}=0$
are parameterized by the restricted cohomology space $H^2_*(\g)$
\cite[Theorem~3.3]{Ho}. If $(\phi,\omega)\in C^2_*(\g)$ is a
restricted 2-cocycle, then the corresponding restricted
one-dimensional central extension $E=\g\oplus\F c$ has Lie bracket and
$[p]$-operation defined by
\begin{align}\label{genonedimext}
  \begin{split}
    [g,h] & =[g,h]_{\g}+\phi(g\wedge h) c;\\
    [g, c] & = 0;\\
    g^{[p]} & = p^{[p]_{\g}}+\omega(g) c;\\
    c^{[p]} & = 0,
  \end{split}
\end{align}
where $[\cdot,\cdot]_\g$ and $\cdot^{[p]_\g}$ denote the Lie bracket
and $[p]$-operation in $\g$, respectively. We can use
(\ref{genonedimext}) to explicitly describe the restricted
one-dimensional central extensions corresponding to the restricted
cocycles in Theorems~\ref{H2forpnot2}, \ref{H2forpcong2mupluslambda0}
and \ref{H2forpcong2mupluslambdanot0}. For the rest of this section,
let $g=\sum \alpha_ie_i$ and $h=\sum\beta_ie_i$ denote two arbitrary
elements of $\g^\lambda(p)$ or $\g^{\mu,\lambda}(p)$.

Just as in \cite{EvFi,EvFim2,EvFim0}, if $E_k=g^\lambda(p)\oplus \F c$
denotes the one-dimensional restricted central extension of
$\g^\lambda(p)$ (or $\g^{\mu,\lambda}(p)$) determined by the
cohomology class of the restricted cocycle $(0,\overline e^k)$, then
(\ref{genonedimext}) gives the bracket and $[p]$-operation in $E_k$:
\begin{align*}
  \begin{split}
    [g,h] & =[g,h]_{\g^\lambda(p)};\\
    [g, c] & = 0;\\
    g^{[p]} & = g^{[p]_{\g^\lambda(p)}}+ \alpha_k^p c;\\
    c^{[p]} & = 0.
  \end{split}
\end{align*}

For the restricted cocycles
$(e^{1,4}, \widetilde{e^{1,4}}), (e^{2,5}, \widetilde{e^{2,5}})$ and
$(\phi_{p+1},\widetilde{\phi_{p+1}})$, we summarize the corresponding
restricted one-dimensional central extensions $E$ in the following
Table (see Remark~\ref{tildemaps}). Everywhere in the Table, we omit
the brackets $[g,c]=0$ and $[p]$-operation $c^{[p]}=0$ for brevity.

\begin{equation*}
  \begin{array}{|l|r l|}
    \hline
    \multicolumn{3}{|c|}{\mbox{\rm {\bf Table 1:} Restricted one-dimensional
    central extensions}}\\
\multicolumn{3}{|c|}{\mbox{\rm with $p\equiv 2\pmod* 3$ and $\mu+\lambda=0$}}\\
    \hline
    \multicolumn{3}{|c|}{p>5}\\
    \hline
 &
   [g,h]=&[g,h]_{\g^{\mu,\lambda}(p)}+(\alpha_1\beta_4-\alpha_4\beta_1) c\\
    (e^{1,4},0)&  g^{[p]} = &g^{[p]_{\g^{\mu,\lambda}(p)}}\\
    \hline
 &
   [g,h]=&[g,h]_{\g^{\mu,\lambda}(p)}+(\alpha_2\beta_5-\alpha_5\beta_2) c\\
    (e^{2,5},0)&  g^{[p]} = &g^{[p]_{\g^{\mu,\lambda}(p)}}\\
    \hline
 &
   [g,h]=&[g,h]_{\g^{\mu,\lambda}(p)}\\
(\phi_{p+1},\widetilde{\phi_{p+1}})& & +\left(\sum_{i=1}^{s(p+1)}a_{i,p+1-i}(\alpha_i\beta_{p+1-i}-\alpha_{p+1-i}\beta_i)\right) c\\
    &  g^{[p]} = &g^{[p]_{\g^{\mu,\lambda}(p)}}+\alpha_1^{p-1}\alpha_2c\\
    \hline
  \end{array}
\end{equation*}
If $p=5$, the extensions $(e^{1,4},0)$ and
$ (\phi_{p+1},\widetilde{\phi_{p+1}})$ are as in Table~1, but the
map $\widetilde{e^{2,5}}\ne 0$ and the $[p]$-operator for the extension
$(e^{2,5},\widetilde{e^{2,5}})$ is given by
\[ g^{[p]}
  =g^{[p]_{\g^{\mu,\lambda}(p)}}+\frac{1}{2}\alpha_1^3\alpha_2^2.\]
The bracket for $(e^{2,5},\widetilde{e^{2,5}})$ is unchanged.

If $\mu+\lambda\ne 0$, then Theorem~\ref{H2forpcong2mupluslambdanot0}
states $(\phi_{p+1},\widetilde{\phi_{p+1}})$ is not a cocycle. The
other two restricted one-dimensional extensions in Table~1 are unchanged.

If $p\not\equiv 2\pmod* 3$, then $p>5$ and
Theorem~\ref{H2forpcong2mupluslambda0} implies
$H^2_*(\g^{\mu,\lambda}(p))=H^2_*(\g^\lambda(p))$ so that the
restricted one-dimensional central extensions are the same as those in
Table 1 with $p>5$ and $\g^{\mu,\lambda}(p)$ replaced with
$\g^\lambda(p)$.

\bibliography{references}{}

\begin{thebibliography}{10}

\bibitem{ChE}
C.~Chevalley and S.~Eilenberg.
\newblock Cohomology theory of {L}ie groups and {L}ie algebras.
\newblock {\em Trans. Amer. Math. Soc.}, 63:84--124, 1948.

\bibitem{EvFi}
T.J. Evans and A.~Fialowski.
\newblock Restricted one-dimensional central extensions of restricted simple
  {L}ie algebras.
\newblock {\em Linear Algebra and Its Applications}, 513C:96--102, 2017.

\bibitem{EvFim2}
T.J. Evans and A.~Fialowski.
\newblock Cohomology of restricted filiform {L}ie algebras
  ${\bf\mathfrak{m}}^\lambda_2(p)$.
\newblock {\em SIGMA}, 15(095):11 pages, 2019.

\bibitem{EvFim0}
T.J. Evans and A.~Fialowski.
\newblock Restricted one-dimensional central extensions of the restricted
  filiform {L}ie algebras ${\bf\mathfrak{m}}^\lambda_0(p)$.
\newblock {\em Linear Algebra and Its Applications}, 595:244--257, March 2019.

\bibitem{EvFiPe}
T.J. Evans, A.~Fialowski, and M.~Penkava.
\newblock Restricted cohomology of modular {W}itt algebras.
\newblock {\em Proc. Amer. Math. Soc.}, 144:1877--1886, 2016.

\bibitem{EvansFuchs2}
T.J. Evans and D.B. Fuchs.
\newblock A complex for the cohomology of restricted {L}ie algebras.
\newblock {\em Journal of Fixed Point Theory and Applications}, 3(1):159--179,
  May 2008.

\bibitem{Feldross}
J.~Feldross.
\newblock On the cohomology of restricted {L}ie algebras.
\newblock {\em Comm. Algebra}, 19(10):2865--2906, 1991.

\bibitem{Fi}
A.~Fialowski.
\newblock On the classification of graded {L}ie algebras with two generators.
\newblock {\em Moscow Univ. Math. Bull.}, 38(2):76--79, 1983.

\bibitem{FuchsBook}
D.B. Fuchs.
\newblock {\em Cohomology of Infinite Dimensional {L}ie Algebras}.
\newblock Consultants Bureau, 1986.

\bibitem{Ho}
G.~Hochschild.
\newblock Cohomology of restricted {L}ie algebras.
\newblock {\em American Journal of Mathematics}, 76(3):555--580, 1954.

\bibitem{Jac}
N.~Jacobson.
\newblock {\em {L}ie Algebras}.
\newblock John Wiley, 1962.

\bibitem{Jantzen2006}
J.~C. Jantzen.
\newblock Restricted {L}ie algebra cohomology.
\newblock In {\em Lect. Notes in Math.}, page 1271. 2006.

\bibitem{Jantzen}
J.C. Jantzen.
\newblock {\em Representations of Algebraic Groups}, volume 107 of {\em
  Mathematical surveys and monographs}.
\newblock American Mathematical Society, 2nd edition, 2003.

\bibitem{Kac}
Victor~G. Kac.
\newblock {\em Infinite-dimensional {L}ie algebras}.
\newblock Cambridge University Press, Cambridge, 1985.

\bibitem{StradeFarnsteiner}
H.~Strade and R.~Farnsteiner.
\newblock {\em Modular {L}ie algebras and their representations}, volume 116.
\newblock Marcel Dekker, 1988.

\end{thebibliography}
\bibliographystyle{plain}

\pagebreak
\appendix

\section{Appendix to Proof of Theorem~\ref{gradedk}}

The examples below illustrate the computations in the proof of
Theorem~3.3 for the prime $p=23$ and values of $k=6t+r$ for
$r=0,1,2,3,4,5$. The examples are listed in the order they appear in parts (i), (ii) and
(iii) of the proof. In part (ii), the example for $r=3$ uses $p=19$ because necessarily
$p\equiv 1\pmod* 3$ in this case.

In each matrix, only the rows for $u=1$ and $u=2$
are shown. The circles indicate entries for which there is exactly
one non-zero entry in the row which shows the corresponding coefficient
$\alpha_{i,k-i}=0$. The shading indicates terms from the
simplified equations from (12). 

{\footnotesize
\paragraph{(i) ($[r]_3=1$) $k=28=6\cdot 4+4$ with $r=4$.}
\[\begin{NiceArray}{l|rrrrrrrrr}
\CodeBefore [create-cell-nodes]
\tikz \node [highlight = (5-3) (5-4)] {} ;
\tikz \node [highlight = (8-6) (8-7)] {} ;
\tikz \node [highlight = (11-9) (11-10)] {} ;
\Body
&e^{5,23} &e^{6,22} &e^{7,21} &e^{8,20} &e^{9,19} &e^{10,18}
    &e^{11,17} &e^{12,16} &e^{13,15} \\
\hline \\
e^{1,4,23}&0 & & & & & & & & \\
e^{1,5,22}&0 &1 & & & & & & & \\
e^{1,6,21}& &-1 &-1 & & & & & & \\
e^{1,7,20}& & &1 &0 & & & & & \\
e^{1,8,19}& & & &0 &1 & & & & \\
e^{1,9,18}& & & & &-1 &-1 & & & \\
e^{1,10,17}& & & & & &1 &0 & & \\
e^{1,11,16}& & & & & & &0 &1 & \\
e^{1,12,15}& & & & & & & &-1 &-1 \\
e^{1,13,14}& & & & & & & & &1 \\
e^{2,3,23}&1 & & & & & & & & \\
e^{2,4,22}& &-1 & & & & & & & \\
e^{2,5,21}&1 & & 0& & & & & & \\
e^{2,6,20}& &0 & &1 & & & & & \\
e^{2,7,19}& & &-1 & &-1 & & & & \\
e^{2,8,18}& & & &1 & &0 & & & \\
e^{2,9,17}& & & & & 0& & 1& & \\
e^{2,10,16}& & & & & &-1 & &-1 & \\
e^{2,11,15}& & & & & & &1 & &0 \\
e^{2,12,14}& & & & & & & &0 & \\
\vdots& & & & & & & & & \\
\CodeAfter
   \tikz \draw (4-3) circle (2.5mm) ;
   \tikz \draw (7-6) circle (2.5mm) ;
   \tikz \draw (10-9) circle (2.5mm) ;
\tikz \draw (13-2) circle (2.5mm) ;
\tikz \draw (16-5) circle (2.5mm) ;
\tikz \draw (19-8) circle (2.5mm) ;
\end{NiceArray}
\]

\paragraph{(i)  ($[r]_3=-1$) $k=29=6\cdot 4+5$ with $r=5$.}

\[\begin{NiceArray}{c|rrrrrrrrr}
\CodeBefore [create-cell-nodes]
\tikz \node [highlight = (6-4) (6-5)] {} ;
\tikz \node [highlight = (9-7) (9-8)] {} ;
\tikz \node [highlight = (20-8) (20-10)] {} ;
\Body
&e^{6,23} &e^{7,22} &e^{8,21} &e^{9,20} &e^{10,19} &e^{11,18}
    &e^{12,17} &e^{13,16} &e^{14,15} \\
\hline \\
e^{1,5,23}&1 & & & & & & & & \\
e^{1,6,22}&0 &-1 & & & & & & & \\
e^{1,7,21}& &-1 &0 & & & & & & \\
e^{1,8,20}& & &1 &1 & & & & & \\
e^{1,9,19}& & & &0 &-1 & & & & \\
e^{1,10,18}& & & & &-1 &0 & & & \\
e^{1,11,17}& & & & & &1 &1 & & \\
e^{1,12,16}& & & & & & &0 &-1 & \\
e^{1,13,15}& & & & & & & &-1 &0 \\
e^{2,4,23}&-1 & & & & & & & & \\
e^{2,5,22}& &0 & & & & & & & \\
e^{2,6,21}&1 & & 1& & & & & & \\
e^{2,7,20}& &0 & &-1 & & & & & \\
e^{2,8,19}& & &-1 & &0 & & & & \\
e^{2,9,18}& & & &1 & &1 & & & \\
e^{2,10,17}& & & & & 0& & -1& & \\
e^{2,11,16}& & & & & &-1 & &0 & \\
e^{2,12,15}& & & & & & &1 & &1 \\
e^{2,13,14}& & & & & & & &0 & 1\\
\vdots& & & & & & & & & \\
\CodeAfter
   \tikz \draw (4-3) circle (2.5mm) ;
   \tikz \draw (7-6) circle (2.5mm) ;
   \tikz \draw (10-9) circle (2.5mm) ;
\tikz \draw (12-2) circle (2.5mm) ;
\tikz \draw (15-5) circle (2.5mm) ;
\tikz \draw (18-8) circle (2.5mm) ;
\end{NiceArray}
\]

\paragraph{(i) ($[r]_3=0$ and $[p]_3=-1$) $k=30=6\cdot 5+0$ with $r=0$.}

\[\begin{NiceArray}{c|rrrrrrrr}
\CodeBefore [create-cell-nodes]
\tikz \node [highlight = (5-3) (5-4)] {} ;
\tikz \node [highlight = (5-4) (6-4)] {} ;
\tikz \node [highlight = (6-4) (6-5)] {} ;
\tikz \node [highlight = (8-6) (8-7)] {} ;
\tikz \node [highlight = (8-7) (9-7)] {} ;
\tikz \node [highlight = (9-7) (9-8)] {} ;
\tikz \node [highlight = (14-2) (14-4)] {} ;
\tikz \node [highlight = (17-5) (17-7)] {} ;
\Body
&e^{7,23} &e^{8,22} &e^{9,21} &e^{10,20} &e^{11,19} &e^{12,18}
    &e^{13,17} &e^{14,16} \\
\hline \\
e^{1,6,23}&-1 & & & & & & &  \\
e^{1,7,22}&0 &0 & & & & & & \\
e^{1,8,21}& &-1 &1 & & & & &  \\
e^{1,9,20}& & &1 &-1 & & & &  \\
e^{1,10,19}& & & &0 &0 & & &  \\
e^{1,11,18}& & & & &-1 &1 & &  \\
e^{1,12,17}& & & & & &1 &-1 &  \\
e^{1,13,16}& & & & & & &0 &0  \\
e^{1,14,15}& & & & & & & &-1  \\
e^{2,5,23}&0 & & & & & & &  \\
e^{2,6,22}& &1 & & & & & &  \\
e^{2,7,21}&1 & & -1& & & & &  \\
e^{2,8,20}& &0 & &0 & & & &  \\
e^{2,9,19}& & &-1 & &1 & & &  \\
e^{2,10,18}& & & &1 & &-1 & &  \\
e^{2,11,17}& & & & & 0& & 0&  \\
e^{2,12,16}& & & & & &-1 & &1  \\
e^{2,13,15}& & & & & & &1 &  \\
\vdots& & & & & & & &  \\
\CodeAfter
   \tikz \draw (3-2) circle (2.5mm) ;
   \tikz \draw (19-7) circle (2.5mm) ;
   \tikz \draw (19-9) circle (2.5mm) ;
\end{NiceArray}
\]

\paragraph{(ii) ($r=1$ so $[p]_3=-1$) $k=25=6\cdot 4+1$.}

\[\begin{NiceArray}{c|rrrrrrrrrrr}
\CodeBefore [create-cell-nodes]
\tikz \node [highlight = (13-2) (13-5)] {} ;
\tikz \node [highlight = (15-2) (15-5)] {} ;
\Body
&e^{2,23} &e^{3,22} &e^{4,21} &e^{5,20} &e^{6,19} &e^{7,18}
    &e^{8,17} &e^{9,16} &e^{10,15} & e^{11,14} & e^{12,13}\\
\hline \\
e^{1,2,22}&0 &1 & & & & & & & & & \\
e^{1,3,21}& &-1 &-1 & & & & & & & & \\
e^{1,4,20}& &  &1 &0 & & & & & & & \\
e^{1,5,19}& & & &0 & 1& & & & & & \\
e^{1,6,18}& & & & &-1 & -1& & & & &\\
e^{1,7,17}& & & & & &1 & 0& & & &\\
e^{1,8,16}& & & & & & &0 & 1& & &\\
e^{1,9,15}& & & & & & & &-1 &-1 & &\\
e^{1,10,14}& & & & & & & & &1 &0 &\\
e^{1,11,13}& & & & & & & & & & 0&1\\
e^{2,3,20}&1 &0 & &1 & & & & & & &\\
e^{2,4,19}&0 & & -1& &-1 & & & & & &\\
e^{2,5,18}&-1 & & & 1& & 0& & & & &\\
e^{2,6,17}&1 & & & &0 & &1 & & & &\\
e^{2,7,16}& 0 & & & & &0 & & -1& & &\\
e^{2,8,15}&-1 & & & & & & 1& &0 & &\\
e^{2,9,14}&1 & & & & & & & 0& &1 &\\
e^{2,10,13}&0 & & & & & & & & -1& &-1\\
e^{2,11,12}& -1 & & & & & & & & &1 &0\\
\vdots& & & & & & & & & & &\\
\end{NiceArray}
\]

\pagebreak

\paragraph{(ii) ($r=3$ so $[p]_3=1$, $p=19$), $k=21=6\cdot 3+3$.}
\[
\]
\[\begin{NiceArray}{c|rrrrrrrrr}
\CodeBefore [create-cell-nodes]
\tikz \node [highlight = (3-2) (3-3)] {} ;
\tikz \node [highlight = (3-3) (4-3)] {} ;
\tikz \node [highlight = (4-3) (4-4)] {} ;
\tikz \node [highlight = (6-5) (6-6)] {} ;
\tikz \node [highlight = (6-6) (7-6)] {} ;
\tikz \node [highlight = (7-6) (7-7)] {} ;
\tikz \node [highlight = (9-8) (9-9)] {} ;
\tikz \node [highlight = (9-9) (10-9)] {} ;
\tikz \node [highlight = (10-9) (10-10)] {} ;
\tikz \node [highlight = (11-2) (11-5)] {} ;
\tikz \node [highlight = (14-2) (14-8)] {} ;
\tikz \node [highlight = (17-2) (17-10)] {} ;
\Body
&e^{2,19} &e^{3,18} &e^{4,17} &e^{5,16} &e^{6,15} &e^{7,14}
    &e^{8,13} &e^{9,12} & e^{10,11} \\
\hline \\
e^{1,2,18}&-1 &1 & & & & & & & \\
e^{1,3,17}& &1 &-1 & & & & & &\\
e^{1,4,16}& & &0 &0 & & & & & \\
e^{1,5,15}& & & &-1 &1 & & & & \\
e^{1,6,14}& & & & &1 &-1 & & & \\
e^{1,7,13}& & & & & &0 & 0& & \\
e^{1,8,12}& & & & & & &-1 &1 & \\
e^{1,9,11}& & & & & & & &1 &-1 \\
e^{2,3,16}&-1 & -1 & & 1& & & & & \\
e^{2,4,15}&1 & &1 & &-1 & & & & \\
e^{2,5,14}&0 & & &0 & &0 & & & \\
e^{2,6,13}& -1& & & &-1 & &1 & & \\
e^{2,7,12}& 1& & & & & 1& & -1& \\
e^{2,8,11}&0 & & & & & & 0& & 0\\
e^{2,9,10}&-1 & & & & & & &-1 &-1 \\
\vdots& & & & & & & & & \\
 % \CodeAfter
%  \OverBrace[shorten,yshift=3pt]{1-2}{18-4}{A_1}
% \OverBrace[shorten,yshift=3pt]{1-5}{18-7}{A_2}
% \OverBrace[shorten,yshift=3pt]{1-8}{18-10}{A_3}
\end{NiceArray}
\]

\paragraph{(iii) ($r=1$) $k=19=6\cdot 3+1$.}

\[\begin{NiceArray}{c|rrrrrrrr|r}
\CodeBefore [create-cell-nodes]
\tikz \node [highlight = (13-3) (13-9)] {} ;
\tikz \node [highlight = (15-3) (15-9)] {} ;
\Body
& e^{1,18} & e^{2,17} &e^{3,16} &e^{4,15} &e^{5,14} &e^{6,13} &e^{7,12} &e^{8,11} &e^{9,10} \\
\hline\\
e^{1,2,16}&1 & 0 & 1 & & & & & & \\
e^{1,3,15}&0 & & -1& -1& & & & & \\
e^{1,4,14}&-1 & & &1 &0 & & & & \\
e^{1,5,13}&1& & & &0 &1 & & &  \\
e^{1,6,12}&0 & & & & & -1&-1 & & \\
e^{1,7,11}& -1& & & & & &1 &0 & \\
e^{1,8,10}&1 & & & & & & &0 &1 \\
e^{2,3,14}& &1 &0 & &1 & & & & \\
e^{2,4,13}& &0 & &-1 & &-1 & & & \\
e^{2,5,12}& &-1 & & & 1& & 0& & \\
e^{2,6,11}& &1 & & & &0 & & 1& \\
e^{2,7,10}& &0 & & & & & -1& &-1 \\
e^{2,8,9}& &-1 & & & & & & 1& 0\\
\vdots& & & & & & & & & \\
\CodeAfter
 \tikz \draw (9-2) circle (2.5mm) ;
\end{NiceArray}
\]

\paragraph{(iii) ($r=2$) $k=20=6\cdot 3+2$.}

\[\begin{NiceArray}{c|rrrrrrrr|r}
& e^{1,19} & e^{2,18} &e^{3,17} &e^{4,16} &e^{5,15} &e^{6,14} &e^{7,13} &e^{8,12} &e^{9,11} \\
\hline\\
e^{1,2,17}&0 & 1 & 1 & & & & & & \\
e^{1,3,16}&-1 & & 0& -1& & & & & \\
e^{1,4,15}&1 & & &-1 &0 & & & & \\
e^{1,5,14}&0& & & &1 &1 & & &  \\
e^{1,6,13}&-1 & & & & & 0&-1 & & \\
e^{1,7,12}& 1& & & & & &-1 &0 & \\
e^{1,8,11}&0 & & & & & & &1 &1 \\
e^{1,9,10}&-1 & & & & & & & &0 \\
e^{2,3,15}& &0 &1 & &1 & & & & \\
e^{2,4,14}& &-1 & &0 & &-1 & & & \\
e^{2,5,13}& &1 & & & -1& & 0& & \\
e^{2,6,12}& &0 & & & &1 & & 1& \\
e^{2,7,11}& &-1 & & & & &0 & &-1 \\
e^{2,8,10}& &1 & & & & & & -1& \\
\vdots& & & & & & & & & \\
\CodeAfter
 \tikz \draw (10-2) circle (2.5mm) ;
\tikz \draw (15-3) circle (2.5mm) ;
\end{NiceArray}
\]

\paragraph{(iii) ($r=4$) $k=22=6\cdot 3+4$.}

\[\begin{NiceArray}{c|rrrrrrrrr|r}
& e^{1,21} & e^{2,20} &e^{3,19} &e^{4,18} &e^{5,17} &e^{6,16} &e^{7,15} &e^{8,14} &e^{9,13}& e^{10,12}\\
\hline\\
e^{1,2,19}&1 & 0 & 1 & & & & & && \\
e^{1,3,18}&0 & & -1& -1& & & & && \\
e^{1,4,17}&-1 & & &1 &0 & & & && \\
e^{1,5,16}&1& & & &0&1 & & &&  \\
e^{1,6,15}&0 & & & & & -1&-1 & && \\
e^{1,7,14}& -1& & & & & &1 &0 && \\
e^{1,8,13}&1 & & & & & & &0 &1& \\
e^{1,9,12}&0 & & & & & & & &-1&-1 \\
e^{1,10,11}&-1 & & & & & & & &&1 \\
e^{2,3,17}& &1 &0 & &1 & & & && \\
e^{2,4,16}& &0 & &-1 & &-1 & & && \\
e^{2,5,15}& &-1 & & & 1& & 0& && \\
e^{2,6,14}& &1 & & & &0 & & 1& &\\
e^{2,7,13}& &0 & & & & &-1 & &-1 &\\
e^{2,8,12}& &-1 & & & & & & 1& &0\\
e^{2,9,11}& &1 & & & & & & & 0&\\
\vdots& & & & & & & & & &\\
\CodeAfter
 \tikz \draw (11-2) circle (2.5mm) ;
\tikz \draw (18-3) circle (2.5mm) ;
\end{NiceArray}
\]

\paragraph{(iii) ($r=5$) $k=23=6\cdot 3+5$.}

\[\begin{NiceArray}{c|rrrrrrrrrr|r}
\CodeBefore [create-cell-nodes]
\tikz \node [highlight = (10-2) (10-11)] {} ;
\tikz \node [highlight = (11-2) (11-11)] {} ;
\Body
& e^{1,22} & e^{2,21} &e^{3,20} &e^{4,19} &e^{5,18} &e^{6,17} &e^{7,16} &e^{8,15} &e^{9,14}& e^{10,13}&e^{11,12}\\
\hline\\
e^{1,2,20}&0 & 1 & 1 & & & & & &&& \\
e^{1,3,19}&-1 & & 0& -1& & & & &&& \\
e^{1,4,18}&1 & & &-1 &0 & & & &&& \\
e^{1,5,17}&0& & & &1&1 & & && & \\
e^{1,6,16}&-1 & & & & & 0&1 & &&& \\
e^{1,7,15}& 1& & & & & &-1 &0 &&& \\
e^{1,8,14}&0 & & & & & & &1 &1&& \\
e^{1,9,13}&-1 & & & & & & & &0&-1& \\
e^{1,10,12}&1 & & & & & & & &&-1&0 \\
e^{2,3,18}& &0 &1 & &1 & & & &&& \\
e^{2,4,17}& &-1 & &0 & &1 & & && &\\
e^{2,5,16}& &1 & & & -1& & 0& &&& \\
e^{2,6,15}& &0 & & & &1 & & 1& &&\\
e^{2,7,14}& &-1 & & & & &0 & &-1 &&\\
e^{2,8,13}& &1 & & & & & & -1& &0&\\
e^{2,9,12}& &0 & & & & & & & 1&&1\\
e^{2,10,11}& &-1 & & & & & & & &0&1\\
\vdots& & & & & & & & & &\\
\CodeAfter
 \tikz \draw (19-3) circle (2.5mm) ;
\end{NiceArray}
\]

\pagebreak

\paragraph{(iii) ($r=0$) $k=24=6\cdot 4+0$.}

\[
\]

\[\begin{NiceArray}{c|rrrrrrrrrr|r}
\CodeBefore [create-cell-nodes]
\tikz \node [highlight = (3-3) (3-4)] {} ;
\tikz \node [highlight = (3-4) (4-4)] {} ;
\tikz \node [highlight = (4-4) (4-5)] {} ;
\tikz \node [highlight = (6-6) (6-7)] {} ;
\tikz \node [highlight = (6-7) (7-7)] {} ;
\tikz \node [highlight = (7-7) (7-8)] {} ;
\tikz \node [highlight = (9-9) (9-10)] {} ;
\tikz \node [highlight = (9-10) (10-10)] {} ;
\tikz \node [highlight = (10-10) (10-11)] {} ;
\tikz \node [highlight = (13-3) (13-6)] {} ;
\tikz \node [highlight = (16-3) (16-9)] {} ;
\tikz \node [highlight = (19-3) (19-10)] {} ;
\Body
& e^{1,23} & e^{2,22} &e^{3,21} &e^{4,20} &e^{5,19} &e^{6,18} &e^{7,17} &e^{8,16} &e^{9,15}& e^{10,14}&e^{11,13}\\
\hline\\
e^{1,2,21}&-1 & -1 & 1 & & & & & &&& \\
e^{1,3,20}&1 & & 1& -1& & & & &&& \\
e^{1,4,19}&0 & & &0 &0 & & & &&& \\
e^{1,5,18}&-1& & & &-1&1 & & && & \\
e^{1,6,17}&1 & & & & & 1&-1 & &&& \\
e^{1,7,16}& 0& & & & & &0 &0 &&& \\
e^{1,8,15}&-1 & & & & & & &-1 &1&& \\
e^{1,9,14}&1 & & & & & & & &1&-1& \\
e^{1,10,13}&0 & & & & & & & &&0&0 \\
e^{1,11,12}&-1 & & & & & & & &&&-1\\
e^{2,3,19}& &-1 &-1 & &1 & & & &&& \\
e^{2,4,18}& &1 & &1 & &-1 & & && &\\
e^{2,5,17}& &0 & & & 0& & 0& &&& \\
e^{2,6,16}& &-1 & & & &-1 & & 1& &&\\
e^{2,7,15}& &1 & & & & &1 & &-1 &&\\
e^{2,8,14}& &0 & & & & & & 0& &0&\\
e^{2,9,13}& &-1& & & & & & & -1&&1\\
e^{2,10,12}& &1 & & & & & & & &1&\\
\vdots& & & & & & & & & &\\
\CodeAfter
 \tikz \draw (12-2) circle (2.5mm) ;
 % \OverBrace[shorten,yshift=3pt]{1-3}{last-5}{A_1}
% \OverBrace[shorten,yshift=3pt]{1-6}{21-8}{A_2}
% \OverBrace[shorten,yshift=3pt]{1-9}{21-11}{A_3}
\end{NiceArray}
\]

\pagebreak

\paragraph{(iii) ($r=3$)  $k=21=6\cdot 3+3$.}

{\normalsize \paragraph{} Here the shaded rows indicate $\alpha_{1,k-1}=-\alpha_{s(k)-1,k-s(k)+1}=\alpha_{2,k-2}$
as in (\ref{theinductionseed}), and the single shaded entries indiate
the terms $\alpha_{i+1,k-i-1}$ and $\alpha_{i+2,k-i-2}$ from
(\ref{induction}).}

\[\begin{NiceArray}{c|rrrrrrrrr|r}
       \CodeBefore [create-cell-nodes]
    \tikz \node [highlight = (10-2) (10-10)] {} ;
    \tikz \node [highlight = (17-3) (17-10)] {} ;
    \tikz \node [highlight = (3-4) (3-4)] {} ;
    \tikz \node [highlight = (4-5) (4-5)] {} ;
    \tikz \node [highlight = (11-6) (11-6)] {} ;
    \tikz \node [highlight = (6-7) (6-7)] {} ;
    \tikz \node [highlight = (7-8) (7-8)] {} ;
    \tikz \node [highlight = (14-9) (14-9)] {} ;
    \tikz \node [highlight = (9-10) (9-10)] {} ;
    \Body
& e^{1,20} & e^{2,19} &e^{3,18} &e^{4,17} &e^{5,16} &e^{6,15} &e^{7,14} &e^{8,13} &e^{9,12}& e^{10,11}\\
\hline\\
e^{1,2,18}&-1 & -1 & 1 & & & & & && \\
e^{1,3,17}&1 & & 1& -1& & & & && \\
e^{1,4,16}&0 & & &0&0 & & & && \\
e^{1,5,15}&-1& & & &-1&1 & & &&  \\
e^{1,6,14}&1 & & & & & 1&-1 & && \\
e^{1,7,13}& 0& & & & & &0 &0 && \\
e^{1,8,12}&-1 & & & & & & &-1 &1& \\
e^{1,9,11}&1 & & & & & & & &1&-1 \\
e^{2,3,16}& &-1 &-1 & &1 & & & && \\
e^{2,4,15}& &1 & &1 & &-1 & & && \\
e^{2,5,14}& &0 & & & 0& & 0& && \\
e^{2,6,13}& &-1 & & & &-1 & & 1& &\\
e^{2,7,12}& &1 & & & & &1 & &-1 &\\
e^{2,8,11}& &0 & & & & & & 0& &0\\
e^{2,9,10}& &-1 & & & & & & & -1&-1\\
\vdots& & & & & & & & & &\\
\end{NiceArray}
\]
}

\end{document}